\documentclass{amsart}
\usepackage{amssymb}
\usepackage{enumerate}
\def\sact{*}

\def\np{\bigskip\noindent}
\def\nl{\smallskip\noindent}

\def\Supp{{\rm Supp}}

\newtheorem{Thm}{Theorem}[section]
\newtheorem{Lm}[Thm]{Lemma}

\newtheorem{lemma}[Thm]{Lemma}
\newtheorem{Prop}[Thm]{Proposition}
\newtheorem{Cor}[Thm]{Corollary}

\theoremstyle{definition}

\newtheorem{Example}[Thm]{Example}

\newtheorem{Alg} [Thm]{Algorithm}

\renewcommand{\phi}{\varphi}

\DeclareMathOperator{\het}{ht}
\def\Z{{\mathbb Z}}

\def\R{{\mathbb R}}

\def\a{\alpha}
\def\b{\beta}
\def\c{\gamma}
\def\s{\sigma}

\def\alp{\alpha}

\author{Arjeh M. Cohen
\& David B. Wales}
\address{Arjeh M. Cohen\\
Department of Mathematics and Computer Science\\
Eindhoven University of Technology\\
POBox 513\\
5600 MB Eindhoven\\
The Netherlands}
\email{A.M.Cohen@tue.nl}
\address{
David B. Wales\\
Mathematics Department\\
Sloan Lab\\
Caltech\\
Pasadena, CA 91125\\
USA}
\email{dbw@its.caltech.edu}

\thanks{This paper was written during
a stay of the first author at Caltech. 
He wants to thank the institute for its hospitality.}

\title{Linearity of Artin groups of finite type}

\begin{document}

\begin{abstract}
Recent results on the linearity of braid groups are
extended in two ways.  We generalize the Lawrence Krammer
representation as well as Krammer's faithfulness proof for this
linear representation to Artin groups of finite type.
\end{abstract}

\maketitle

\section{Introduction}\label{intro}
Recently, both Bigelow \cite{Bigelow} and Krammer \cite{Krammer}
proved that the braid groups are linear. The braid group on $n+1$
braids is the Artin group of type $A_n$. This paper extends the result
to all Artin groups whose types are finite, that is, belong to finite
Coxeter groups.

\begin{Thm}
\label{mainres}\label{mainthm}
Every Artin group of finite type is linear.
\end{Thm}

\np Linearity of a group means that it has a faithful linear representation.  A
standard argument reduces the proof to Artin groups whose types are
finite and irreducible.  We focus on the Artin groups of type $A$,
$D$, $E$.  Since the other Artin groups of finite and irreducible type
can be embedded in these (cf.~\cite{Crisp}), it will suffice for a
proof of Theorem \ref{mainres} to exhibit a faithful representation
for each of the groups $A$, $D$, $E$. The theorem below provides more
information about the representation found. 

Throughout the paper,
we fix a Coxeter matrix $M$ of dimension $n$, and
denote by $B$ the Artin group of type $M$.  This means that
$B$ is the group generated by
$n$ elements $s_1,\ldots,s_n$ subject to the relations
\begin{eqnarray}\label{Artinrels}
\underbrace{s_is_js_i\cdots}_{\mbox{length }M_{ij}} &=&
\underbrace{s_js_is_j\cdots}_{\mbox{length }M_{ij}}
\end{eqnarray}
for $1\le i< j\le n$.  The Coxeter system of type $M$ is denoted by
$(W,R)$ with $R$ consisting of the images $r_i$ of $s_i$
$(i=1,\ldots,n)$ under the natural homomorphism from $B$ to $W$. 
We use the standard facts and some terminology of root
systems as treated for example in \cite{Bourb}.  We
shall be working solely with Artin groups of finite type, so $W$ is
assumed finite, and $W$ has a finite root system $\Phi$ in $\R^n$.  We
shall denote by $\a_1,\ldots,\a_n$ the fundamental roots,
corresponding to the reflections $r_1,\ldots,r_n$,
respectively, and by $\Phi^+$ the set of positive roots:
$$\Phi^+ = \Phi\cap\bigoplus_{1\le i\le n}\R_{\ge0}\a_i.$$
Then $\Phi$ is the
disjoint union of $\Phi^+$ and $\Phi^-=-\Phi^+$.
If $M=A_n$
$(n\ge1)$, $D_n$ $(n\ge4)$, $E_6$, $E_7$, or $E_8$, we say that $B$ is
of type $A$, $D$, $E$.

The coefficients of our representation will be taken in the
ring $\Z[r,t, r^{-1},t^{-1}]$, and we write $V$ for the
free module over that ring with generators $x_\b$ indexed by $\b\in\Phi^+$.

\begin{Thm} \label{adethm}
Let $B$ be an Artin group of type $A$, $D$, $E$.
Then, for each $k\in\{1,\ldots,n\}$
and each $\b\in\Phi^+$, there is a polynomial $T_{k,\b}$
in $\Z[r]$
such that
the following map on the generators of $B$ determines a representation of $B$
on $V$.
$$s_k\mapsto \sigma_k  =  \tau_k+tT_{k},$$
where $\tau_k$ is determined by
\begin{eqnarray*}
\tau_k(x_\b)&=&
\left\{
\begin{tabular}{lr}
$0$ & \mbox{ if }\quad $\phantom{-}(\a_k,\b)=2$\phantom{,}\\
$r  x_{\b-\a_k}$ & \mbox{ if }\quad $\phantom{-}(\a_k,\b)=1$\phantom{,}\\
$ x_{\b}$ & \mbox{ if }\quad $\phantom{-}(\a_k,\b)=0$\phantom{,}\\
$ (1-r^2)  x_{\b} + r x_{\b+\a_k}$ &\phantom{-}
 \mbox{ if }\quad $(\a_k,\b)=-1$,\\
\end{tabular}
\right.
\end{eqnarray*}
and $T_k$ is the linear map on $V$ determined by $T_kx_\b = T_{k,\b}x_{\a_k}$
on the generators of $V$.
If $r$ is specialized to a real number $r_0$, $0<r_0<1$,
in $V\otimes \R$,
we obtain a faithful representation of $B$ on the
resulting free $\R[t,t^{-1}]$-module
$V_1$ with basis $x_\b$ $(\b\in\Phi^+)$.
\end{Thm}

\np The proof of this theorem is based on Krammer's methods. More
specifically, we generalize the Lawrence Krammer representation
\cite{Krammer,Lawrence} as well as Krammer's faithfulness results from
braid groups to Artin groups corresponding to a spherical root system
with a single root length.  The difficulties in the proof come in in
the proper definition of the $T_{k,\b}$.  These are determined by
Algorithm~\ref{alg}.  The algorithm has been implemented in the
computer algebra package Maple and has been used to construct the
representations of $B$ for all $M$ of type $A$, $D$, $E$, and rank
at most 10.

\np In Section \ref{basics} we recall basic and useful properties of
Artin groups, including generalizations of some of Krammer's results
on braid groups.  Section \ref{rep} introduces the representation
referred to in Theorem \ref{adethm}.  Section \ref{krammer} presents a
version of Krammer's linearity proof generalized to Artin groups of arbitrary
types and applies it to the representation of the preceding section.
The paper finishes with a few remarks on alternative proofs
in Section \ref{epilog}.

\section{Basic properties of Artin groups}
\label{basics}

We maintain the notation of the introduction.
The Coxeter group $W$ is assumed to be finite of type
$A$, $D$, or $E$.

The submonoid of $B$ generated by $s_1,\ldots,s_n$ is denoted by
$B^+$.  
By $\le$ we denote the partial order on $B^+$ given by
$x\le y\Leftrightarrow  y \in xB^+$.
The length function on $W$ with respect to $R$, as well as the
length function on $B$ with respect to $\{s_1,\ldots,s_n\}$, is
denoted by $l$.

\np
\begin{Prop}\label{Artinmonprops}
The Artin monoid $B^+$
satisfies the following properties.
\begin{enumerate}[(i)]
\item
The relations (\ref{Artinrels}) form a presentation for
$B^+$ as a monoid generated by $s_1,\ldots,s_n$.
\item
$B=(B^+)^{-1}B^+$. Consequently, 
if $\rho $ is a faithful
linear representation of the monoid $B^+$
such that $\rho(s_i)$ is invertible,
then $\rho$ extends uniquely to a faithful linear representation of $B$.
\item For $x,y,z\in B^+$ we have
$zx\le zy \Leftrightarrow x\le y$.
\item There is a uniquely determined map $b: W \to B^+$ satisfying
$b(uv) = b(u)b(v)$ whenever $u,v\in W$ with $l(uv)=l(u)+l(v)$.
It is injective and satisfies $l(b(u)) = l(u)$.
Write $\Omega = b(W)\subset B^+$.
\item There is a uniquely determined map
$L: B^+\to W$ such that,
for each $x\in B^+$, its image $L(x)$ is the longest element $w$ of $W$
with the property that $b(w)\le x$.
\item For $x$, $y\in B^+$, we have $L(xy) =L(xL(y))$.
\item The map
$y\mapsto L(xy)$ $(x\in B^+)$ defines an action 
of $B^+$ on $\Omega$.
\end{enumerate}
\end{Prop}

\np
\begin{proof}
(i), the first part of (ii), (iii), and (iv) go back to \cite{Deligne};
they are also stated in \cite{Crisp,Charney}.
The second part of (ii) is a direct consequence of the first part
(observed in \cite{Krammer}).

\nl(v). See \cite{Charney}.

\nl (vi). This is Corollary 1.23 of \cite{Deligne} (cf.\ 
Lemma 2.4 of \cite{Charney}).

\nl (vii). As observed in \cite{Krammer}, this is immediate from (vi).
\end{proof}

\np
A subset $A$ of $\Phi^+$ is called {\em closed}
when 
$$\a,\b\in A,\ \ \a+\b\in\Phi^+\Rightarrow \a+\b\in A.$$
By $\mathcal C$ we denote the collection of all closed subsets of $\Phi^+$.  
For $w\in W$, set
$$\Phi_w = \{\a\in\Phi^+\mid w^{-1}\alp\in\Phi^-\}.$$
Let $\mathcal D$ be the collection of all $\Phi_w$ $(w\in W)$.  

On $W$ we have a partial order given
by
\begin{eqnarray}\label{w=vu}
v\le w &\Leftrightarrow &\exists_{u\in W}\ w=vu\ \mbox{ and }
l(w) = l(u) + l(v).
\end{eqnarray}
On $\mathcal C$, we consider the partial order by inclusion.

\begin{lemma}\label{Phiprops}
The members of $\mathcal D $ have the following properties.
\begin{enumerate}[(i)]
\item If $v$, $u$, $w$ are as in (\ref{w=vu}), then
$\Phi_w = \Phi_v\cup v(\Phi_u)$.
\item The size of $\Phi_w$ equals $l(w)$.
\item The members of $\mathcal D$ are closed.  
\item If $A$ is a closed subset of $\Phi^+$, then there is a unique maximal
subset $A'$ of $A$ of the form $\Phi_w$ with $w\in W$.
\item There is an isomorphism of partially ordered
sets $(W,\le)\to ({\mathcal D},\subseteq)$  given by $w\mapsto \Phi_w$.
\end{enumerate}
\end{lemma}

\np
\begin{proof}
\nl(i) and (ii). See \cite{Bourb}.

\nl(iii). For $\b,\c\in \Phi_w$,
we have $w^{-1}(\b+\c)\in(\Phi^-+\Phi^-)\cap \Phi \subseteq \Phi^-$.

\nl(iv).
For $x,y\in W$ we have $x\le xy$ if and only if $l(xy) = l(x)+l(y)$,
if and only if
$\Phi_{xy} = \Phi_x\cup x(\Phi_y)$, if and only if $\Phi_x\subseteq \Phi_{xy}$.

Suppose now that there is no largest member of $\mathcal D$ contained in
$A$. Then, by (i), there are $u\in W$ and $i,j\in \{1,\ldots,n\}$ with
$u\le ur_i$ and $u\le ur_j$ for which $\Phi_{ur_i} \subseteq A$ and
$\Phi_{ur_j} \subseteq A$ such that no member of $\mathcal D$ containing
$\Phi_{ur_i}\cup\Phi_{ur_j}$ is a subset of $A$.  Then $u\a_i $ and
$u\a_j$ are in $A$, and, by \cite{Bourb}, $u\le uw_{ij}$, where
$w_{ij}$ is the longest element of the subgroup of $W$ generated by
$r_i$ and $r_j$.  This is $r_ir_j$ if $r_i$ and $r_j$ commute and
$r_jr_ir_j$ if they do not.  So, in the former case
$\Phi_{uw_{ij}}\setminus\Phi_w$ consists of $u\a_i$ and $u\a_j$, and
in the latter case of $u\a_i$, $u\a_j$, and $u\a_i+u\a_j$, which
belongs to $A$ as $A$ is closed.  This means $\Phi_{uw_{ij}} \subseteq
A$, a contradiction with $\Phi_{ur_{i}}\cup \Phi_{ur_{j}} \subseteq\Phi_{uw_{ij}}$.

\nl(v). We have $l(w) = |\Phi_w|$ so we can work by induction on
$l(w)$. Clearly, $\Phi_1$ is the empty set, so assume $l(w)>1$.
Then there is a fundamental reflection $r_i$ such that
$l(r_iw)<l(w)$. Now
$w^{-1}(\a_i)\in\Phi^-$ so, if $\a_i\in \Phi_w = \Phi_v$,
then also $v^{-1}(\a_i)\in\Phi^-$, so
$l(r_iv)<l(v)$. Consequently, $\Phi_{r_iw}=\Phi_{r_iv}$,
so by induction $r_iw=r_iv$, establishing $w=v$.
\end{proof}

\np 
For $A$ a closed subset of $\Phi^+$, write $g(A) = x$ for
$x\in\Omega$ such that $\Phi_{b^{-1}x}$ is the maximal subset of $A$
belonging to $\mathcal D$.  In view of Proposition
\ref{Artinmonprops} and Lemma \ref{Phiprops}, the map $g: {\mathcal
C}\to \Omega$ is well defined.

In the next section we define the linear representation for $B$ of
type $A$, $D$, $E$.  In the subsequent section, we use this
representation to define an action of $B^+$ on $\mathcal C$ that makes the
map $g$ equivariant with the action on $\Omega$ of Proposition
\ref{Artinmonprops} (vii).

\section{The representation for types $A$, $D$, $E$}
\label{rep}

In this section, we continue to assume that
the type $M$ is $A$, $D$, or $E$.
This has the consequences that $\Phi^+$ is finite and that $M$ can be
viewed as a graph on $\{1,\ldots,n\}$ with adjacency $k\sim l$
given by $(r_kr_l)^3=1$, or, equivalently, $(\a_k,\a_l) =
-1$. Nonadjacency of $k$ and $l$ corresponds to $(r_kr_l)^2 = 1$ and to
$(\a_k,\a_l) = 0$. 

We shall first describe the `$t=0$ part' of the linear
representation of the Artin monoid $B^+$.
Recall the $\Z[r^{\pm1},t^{\pm1}]$-module $V$
and the linear transformations $\tau_k$ introduced
in Theorem \ref{adethm}.
Denote by $V_0$ the free $\Z[r]$-module with generators
$x_\b$ $(\b\in \Phi^+)$. Thus,
$V_0$ is contained in $V$, and $V$ is obtained from $V_0$ by extending
scalars to $\Z[r^{\pm1},t^{\pm1}]$.

\begin{Lm}\label{monmorphismlm}
There is a monoid homomorphism $ B^+ \to {\rm End}(V_0)$
determined by $s_i\mapsto\tau_i$ $(i=1,\ldots,n)$.
\end{Lm}

\begin{proof}
We must show that, if $\a_i$ and $\a_j$ are not adjacent, then $\tau_i
\tau_j=\tau_j \tau_i$ and, if they are adjacent, then
$\tau_i\tau_j\tau_i=\tau_j\tau_i\tau_j$.  We evaluate the expressions
on each $x_\b$ and show they are equal.  We begin with the case in
which $\b=\a_i$ or $\a_j$.  To be specific, let $\b=\a_i$.  Suppose
first that $i$ and $j$ are not adjacent.  Then $\tau_ix_{a_i}=0$ and
$\tau_jx_{\a_i}=x_{\a_i}$.  Now $\tau_j\tau_ix_{\a_i}=0$,
$\tau_i\tau_jx_{\a_i}=\tau_ix_{\a_i}=0$ and the result holds.
Suppose next that $i$ and $j$ are adjacent.  Then
$\tau_ix_{a_i}=\tau_jx_{\a_j}=0$ and
$\tau_jx_{\a_i}=(1-r^2)x_{\a_i}+rx_{\a_i+\a_j}$.  Now
$$\tau_i\tau_j\tau_ix_{\a_i}=\tau_i\tau_j(0)=0$$ and 
\begin{eqnarray*}
\tau_j\tau_i\tau_jx_{\a_i} &=&  \tau_j\tau_i((1-r^2)x_{\a_i}+rx_{\a_i+\a_j}) 
= \tau_j(0+r^2x_{\a_i+a_j-\a_i})  \\
 &=&r^2 \tau_jx_{\a_j} 
 = 0.
\end{eqnarray*}

We now divide the verifications into the various cases depending on
the inner products $(\a_i,\b)$ and $(\a_j,\b)$.  The table below
describes the images of the vectors $x_\b$ under $\tau_i$ and
$\tau_j$.

\begin{center}
\begin{tabular}{|r|r|c|c|}
\hline
 $(\a_i,\beta)$&$(\a_j,\beta)$&$\tau_ix_\b$&$\tau_jx_\b$\\
\hline
 $0$\ \ &$0$\ \ &$x_\b$&$x_\b$ \\
 $1$\ \ &$1$\ \ &$r x_{\b-\a_i}$&$r  x_{\b-a_j}$  \\
 $1$\ \ &$-1$\ \ &$r x_{\b-\a_i}$&$(1-r^2)x_{\b} + r x_{\b+\a_j} $ \\
 $1$\ \ &$0$\ \ &$r x_{\b-\a_i}$&$x_\b$   \\
 $-1$\ \ &$0$\ \ &$(1-r^2)x_{\b} + r x_{\b+\a_i}$&$x_{\b}$ \\
 $-1$\ \ &$-1$\ \ &$(1-r^2)x_{\b} + r x_{\b+\a_i}$&$   (1-r^2)x_{\b} + r x_{\b+\a_j} $    \\
\hline
\end{tabular}
\end{center}

\np
First assume that $(\a_i,\a_j)=0$.  The computations 
verifying $\tau_i\tau_j=\tau_j\tau_i$ are
straightforward.  We summarize the results in the following table.

\medskip
\begin{center}
\begin{tabular}{|r|r|c|}
\hline
 $(\a_i,\beta)$&$(\a_j,\beta)$&$\tau_i\tau_j x_\b=\tau_j\tau_ix_\b$\\
\hline
 $0$\ \ &$0$\ \ &$x_\b$\\
 $1$\ \ &$1$\ \ & $r^2x_{\b-\a_i-\a_j}$  \\
 $1$\ \ &$-1$\ \ &$r(1-r^2)x_{\b-\a_i}+r^2x_{\b+\a_i-\a_j}$ \\
 $1$\ \ &$0$\ \ &$rx_{\b-\a_i}$   \\
 $-1$\ \ &$0$\ \ &$(1-r^2)x_{\b}+rx_{\b+\a_i}$ \\
 $-1$\ \ &$-1$\ \ &$(1-r^2)^2x_\b+r(1-r^2)(x_{\b+\a_i}+x_{\b+a_j})+r^2x_{\b+\a_i+\a_j}$    \\
\hline
\end{tabular}
\end{center}

\nl
We demonstrate how to derive these expressions by checking the third line.  
\begin{eqnarray*}
\tau_i\tau_jx_{\b} &=& \tau_i((1-r^2)x_\b+rx_{\b +\a_j})    
= (1-r^2)rx_{\b-\a_i}+r^2x_{\b+\a_i-\a_j}. 
\end{eqnarray*}
In the other order,
\begin{eqnarray*}
\tau_j\tau_ix_{\b} &=& \tau_j(rx_{\b-\a_i})   
= (1-r^2)rx_{\b-\a_i}+r^2x_{\b+\a_i-\a_j}. 
\end{eqnarray*}
Recall $(\a_i,\a_j)=0$ and so $(\b-\a_i,\a_j)=(\b,\a_j)=-1$.  

Suppose then that $(\a_i,\a_j)=-1.$  The same situation occurs except the 
computations are sometimes longer and one case does not occur.  This is 
the case where $(\a_i,\b)=(\a_j,\b)=-1$.  For then
$\b+\a_i$ is also a root, and $(\b+\a_i,\a_j)=-1-1=-2$.  This means 
$\b+\a_i=-\a_j$ and $\b$ is not a positive root.  
The table is as follows.
\medskip

\begin{center}
\begin{tabular}{|r|r|c|}
\hline
 $(\a_i,\beta) $&$(\a_j,\beta)$&$\tau_i\tau_j \tau_ix_\b=\tau_j\tau_i\tau_jx_\b$\\
\hline
 $0$\ \ &$0$\ \ &$x_\b$\\
 $1$\ \ &$1$\ \ & $0$  \\
 $1$\ \ &$-1$\ \ &$r(1-r^2)x_{\b-\a_i}+r^2x_\b$ \\
 $1$\ \ &$0$\ \ &$r^2x_{\b-\a_i-\a_j}$  \\
 $-1$\ \ &$0$\ \ &$(1-r^2)x_{\b}+(1-r^2)r x_{\b+\a_i}+r^2x_{\b+\a_i+\a_j}$ \\
 $-1$\ \ &$-1$\ \ &{\rm does not occur}   \\
\hline
\end{tabular}
\end{center}

\nl
As above, these calculations are routine.  
Note that, in the second line, $\b=\a_i+\a_j$.
We do the second from the last case
in detail.  Here, $(\a_i,\b)=-1$ and $(\a_j,\b)=0$.
\begin{eqnarray*}
\tau_j\tau_i\tau_jx_{\b} &=& \tau_j \tau_ix_\b  
= \tau_j((1-r^2)x_\b+rx_{\b+\a_i})   \\ 
  &=&(1-r^2)x_\b+r(1-r^2)x_{\b+\a_i}+r^2x_{\b+\a_i+\a_j}  .
\end{eqnarray*}
In the other order,
\begin{eqnarray*}
\tau_i\tau_j\tau_ix_{\b} &=& \tau_i\tau_j((1-r^2)x_\b+rx_{\b+\a_i})  \\
 &=& \tau_i((1-r^2)x_\b+r(1-r^2)x_{\b+\a_i}+rrx_{\b+\a_i+\a_j}) \\
&=&   (1-r^2)^2x_\b+(1-r^2)rx_{\b+\a_i}+r^2(1-r^2)x_\b+r^2x_{\b+\a_i+\a_j} \\
 &=& (1-r^2)x_\b+r(1-r^2)x_{\b+\a_i}+r^2x_{\b+\a_i+\a_j} .
\end{eqnarray*}
\end{proof}

\np
We next study the possibilities for the parameters $T_{k,\b}$ occurring in
Theorem \ref{adethm}.
Recall that there we defined
$\s_k=\tau_k+tT_k$, where
$T_kx_\b=T_{k,\b}x_{\a_k}$.

\begin{Prop}\label{finalTiRelations}
Set $T_{i,\a_i}=r^4$ for all $i\in\{1,\ldots,n\}$.
For $\sigma_i\mapsto \tau_i+tT_i$ to define a linear representation
of the group $B$ on $V$,
it is necessary and sufficient that the equations in Table
\ref{algtable} are satisfied
for each $k,l=1,\ldots,n$ and each $\beta\in\Phi^+$.
\end{Prop}

\np
\begin{proof}
The $\s_k$ should satisfy the relations
(\ref{Artinrels}) for $s_k$.
Substituting $\tau_k+tT_k$ for $s_k$,
we find relations for the coefficients
of $t^i$ with $i=0,1,2,3$.
The constant part involves only the $\tau_k$. It follows from
Lemma \ref{monmorphismlm} that these equations are satisfied.

The coefficients of $t$ lead to
\begin{eqnarray}\label{genTrepeqcom}
T_k\tau_l = T_k\ \  \mbox{ and }\ \ \ 
T_l\tau_k = T_l \ \  &\mbox{if}&\ \ 
(\a_k,\a_l) = 0,\\
\tau_lT_k \tau_l  + T_l\tau_k \tau_l
+\tau_l\tau_kT_l
=\tau_kT_l \tau_k
 + T_k\tau_l \tau_k  + \tau_k\tau_lT_k
 \ \  &\label{genTrepeqadj}\mbox{if}&\ \ 
(\a_k,\a_l) = -1.
\end{eqnarray}

We focus on the consequences of these equations for the $T_{k,\b}$.
Consider the case where $(\a_k,\a_l)=0$. Then
\begin{center}
\begin{tabular}{|r|c|c|}
\hline
$(\a_k,\beta)$&$\tau_kx_\b$&$T_l\tau_k = T_l$\\
\hline
$0$\ \ &$x_\b$&no condition\\
$1$\ \ &$r x_{\b-\a_k}$&$rT_{l,\b-\a_k}=T_{l,\b}$\\
$-1$\ \ &$(1-r^2)x_{\b} + r 
x_{\b+\a_k}$&$(1-r^2)T_{l,\b} + r T_{l,\b+\a_k}=T_{l,\b}$\\
\hline
\end{tabular}
\end{center}

\nl
Both equations say the same, namely,
\begin{eqnarray}\label{commlaweq}
T_{l,\b} &=& rT_{l,\b-\a_k}\ \ \mbox{if}\ (\a_k,\b)=1
\ \mbox{and}\ (\a_k,\a_l)=0.
\end{eqnarray}

\np
Next, we assume $(\a_k,\a_l)=-1$.
A practical rule is
\begin{eqnarray*}\label{alp2alpteq}
\tau_k\tau_lx_{\alp_k}
&=&
\tau_k((1-r^2)x_{\alp_k}+rx_{\alp_k+\alp_l}) = 
r^2x_{\alp_l}.
\end{eqnarray*}
We distinguish cases according to the values of
$(\a_k,\b)$ and $(\a_l,\b)$. Since each inner product, for distinct roots
is one of $1$, $0$, $-1$, there are six cases to consider up to symmetry
(interchange of $k$ and $l$).
However, the case 
$(\a_k,\b)=(\a_l,\b)= -1$ does not occur.
For then $(\s_k\b,\a_l)=-2$, a contradiction with the fact that 
both $\s_k\b$ and $\a_l$ are positive roots.

For the sake of brevity, let us denote the images of the left hand side and
the right hand side of (\ref{genTrepeqadj}) 
on $x_\b$ by
$LHS$ and $RHS$, respectively. 

\nl
Case $(\a_k,\b)=(\a_l,\b)= 1$. Then
$(\s_k\b,\a_l)=(\b-\a_k,\a_l)=2$, so
$\b = \a_k+\a_l$.
Now
\begin{eqnarray*}
LHS& =&
\tau_kT_l \tau_kx_\b  + T_k\tau_l \tau_kx_\b
+\tau_k\tau_lT_kx_\b\\
& =&
r \tau_kT_l x_{\a_l}  + rT_k\tau_l x_{\a_l}
+T_{k,\a_k+\a_l}\tau_k\tau_lx_{\a_k}\\
& =&
T_{l,\a_l}r\tau_kx_{\a_l}  + 
T_{k,\a_k+\a_l} r^2x_{\a_l}\\
& =&
T_{l,\a_l}r(1-r^2)x_{\a_l}  + T_{l,\a_l}r^2 x_{\a_k+\a_l}
+T_{k,\a_k+\a_l} r^2x_{\a_l}\\
& =&
(T_{l,\a_l}r(1-r^2) +T_{k,\a_k+\a_l} r^2) x_{\a_l} 
 + T_{l,\a_l}r^2 x_{\a_k+\a_l} .
\end{eqnarray*}
Comparison with the same expression but then
$l$ and $k$ interchanged yields $RHS$. This leads to the following
two equations.
\begin{eqnarray}\label{11eq1}
T_{k,\a_k+\a_l}  &=& T_{l,\a_l}(r-r^{-1}),\\
T_{k,\a_k}  &=& T_{l,\a_l} .\nonumber
\end{eqnarray}
The second one, and homogeneity of the presentation relations,
allow us to scale the $T_i$ so that 
\begin{eqnarray}\label{11eq2}
T_{i,\a_i} &=& r^4.
\end{eqnarray}

\bigskip\noindent
Case $(\a_k,\b)=(\a_l,\b)= 0$. This gives
\begin{eqnarray*}
LHS& =&
\tau_kT_l x_\b + T_k x_\b
+\tau_k\tau_lT_kx_\b\\
&=&T_{l,\b}\tau_k x_{\a_l} + T_{k,\b} x_{\a_k}
+T_{k,\b}\tau_k\tau_lx_{\a_k}\\
&=&T_{l,\b}(1-r^2) x_{\a_l} +
T_{l,\b}r x_{\a_l+\a_k} + T_{k,\b} x_{\a_k}
+T_{k,\b}r^2 x_{\a_l}\\
&=&(T_{k,\b}r^2 +T_{l,\b}(1-r^2) )x_{\a_l} +
T_{l,\b}r x_{\a_l+\a_k} + T_{k,\b} x_{\a_k}
\end{eqnarray*}
and
$RHS$ can be obtained from the above
by interchanging the indices $k$ and $l$.
Comparison of each of the coefficients of
$x_{\a_k}$, $x_{\a_l+\a_k}$, $x_{\a_l}$ gives

\begin{eqnarray}\label{00eq}
T_{k,\b} & = & T_{l,\b} .
\end{eqnarray}

\np
Since the other cases come down to similar computations, we only list
the results.

\bigskip\noindent
Case $(\a_k,\b)=0$, $(\a_l,\b)= -1$.
%
%
%
Then
\begin{eqnarray}\label{0-1eq1}
T_{k,\b+\a_l}
&=&
r^{-1}T_{l,\b}
-T_{k,\b}(r^{-1}-1),
\\
T_{l,\b+\a_l+\a_k}
&=&
T_{k,\b} - r^{-1}(1-r^2)T_{l,\b+\a_l}.\label{0-1eq2}
\end{eqnarray}

\bigskip\noindent
Case $(\a_k,\b)=0$, $(\a_l,\b)= 1$.
%
%
%
%
%
%
Here
\begin{eqnarray}\label{01eq1}
T_{k,\b}
&=&
T_{l,\b-\a_l-\a_k}
+
(r^2-1)r^{-1}T_{k,\b-\a_l},
\\
T_{l,\b}
&=&
rT_{k,\b-\a_l}.\label{01eq2}
\end{eqnarray}

\bigskip\noindent
Case $(\a_k,\b)=1$, $(\a_l,\b)= -1$.
Now
%
%
%
%
%
\begin{eqnarray}\label{1-1eq1}
T_{l,\b}
&=&
r^{-1}T_{k,\b-\a_k}-
r^{-1}(1-r^2)T_{l,\b-\a_k},
\\
T_{l,\b+\a_l}
&=&
T_{k,\b}r,\label{1-1eq2}
\\
T_{k,\b+\a_l}
&=&
T_{l,\b-\a_k} - (1-r^2)r^{-1}
T_{k,\b}.
\label{1-1eq3}
\end{eqnarray}

\np
\begin{center}
\begin{table}
\caption{Equations for $T_{k,\b}$}
\label{algtable}
\begin{tabular}{|l|c|c|}
\hline
$T_{k,\beta}$&condition&reference\\
\hline
$0$& $\b=\a_l$ and $k\ne l$&(\ref{knotleq})  \\
\hline
$r^4$&$\b=\a_k$&(\ref{11eq2})  \\
\hline
$r^5-r^3$&$\b=\a_k+\a_l$&(\ref{11eq1}) \\
\hline
$rT_{k,\b-\a_l}$&$(\a_l,\b)=1$ and $(\a_k,\a_l)=0$&(\ref{commlaweq})\\
\hline
$T_{l,\b-\a_k-\a_l}+(r-r^{-1})T_{k,\b-\a_l}
$&$(\a_k,\b)=0$ and $(\a_l,\b)=1$&(\ref{01eq1})
\\
& and $(\a_k,\a_l)=-1$&\\
\hline
$r^{-1}T_{l,\b-\a_l}+
(r-r^{-1})T_{k,\b-\a_l}
$&$(\a_k,\b)=-1 $ and $(\a_l,\b)=1$& (\ref{1-1eq1}) 
\\
& and $(\a_k,\a_l)=-1$&\\
\hline
$rT_{l,\b-\a_k}$&$(\a_k,\b)=1$ and $(\a_l,\b)=0$&(\ref{01eq2}) 
\\
& and $(\a_k,\a_l)=-1$&\\
\hline
\end{tabular}
\end{table}
\end{center}

\np
We see that, in order to be a representation,
the $T_{i,\b}$ have to satisfy the equations 
(\ref{commlaweq})--(\ref{1-1eq3}).
But this system is superfluous in that, when the root in the index of
the left hand side of 
(\ref{0-1eq1}) is set to $\gamma$, we obtain 
(\ref{1-1eq1}) for $\gamma$ instead of $\beta$;
similarly,  (\ref{0-1eq2}) is equivalent to
(\ref{01eq1}), while
(\ref{1-1eq2}) is equivalent to
(\ref{01eq2})
and (\ref{1-1eq3}) is equivalent to
(\ref{01eq1}).

We also contend that the equations in (\ref{00eq}) are consequences of the
other relations from Table \ref{algtable}.
The equation says that $T_{k,\beta} = T_{l,\b}$ whenever
$(\a_k,\b)=(\a_l,\b)=0$ and $k\sim l$. We prove this by
induction on the height of $\b$.  The initial cases are direct from
(\ref{11eq1}) and (\ref{11eq2}).  Suppose that $\beta\in\Phi^+$
and $k,l\in\{1,\ldots,n\}$ satisfy the inner product conditions.
There exists $m\in\{1,\ldots,n\}$ such that $(\a_m,\b) = 1$.  
If $(\a_m,\a_k) = (\a_m,\a_l)=0$, then (\ref{commlaweq})
applies to both sides, giving $T_{k,\b} = rT_{k,\b-\a_m}=
rT_{l,\b-\a_m} = T_{l,\b} $, where the middle step uses the induction
hypothesis.

Therefore, interchanging $k$ and $l$ if necessary, we may assume that
$(\a_l,\a_m) = -1$, whence $(\a_k,\a_m) = 0$ (as the Dynkin diagram
contains no triangles).
Now $\delta=\b-\a_m-\a_l\in\Phi^+$ and
$(\a_k,\delta)=1$, so (\ref{commlaweq}) gives
$T_{m,\delta} = rT_{m,\delta-\a_k}$, which, by induction on height, is equal to
$rT_{l,\delta-\a_k}$ (as $(\a_l,\delta-\a_k) = (\a_m,\delta-\a_k) = 0$).
Consequently,
\begin{eqnarray*}
T_{k,\b} &=& rT_{k,\b-\a_m}  \qquad \hfill\phantom{mucchhhhhhhhi}\,
 \mbox{ by (\ref{commlaweq})}\\
 &=& rT_{l,\delta-\a_k} +(r^2-1)T_{k,\delta} \qquad \hfill
\phantom{mx}\, \mbox{ by (\ref{01eq1})}\\
 &=& T_{m,\delta} +(r^2-1)T_{k,\delta} \qquad \hfill
\phantom{mucchi} \mbox{ by the above }\\
 &=& T_{m,\delta} +(r-r^{-1})T_{l,\b-\a_m} \qquad \hfill
\phantom{m} \mbox{ by (\ref{01eq2})}\\
 &=& T_{l,\b} \qquad \hfill \phantom{mxucchmucccchmorre}
\mbox{ by (\ref{01eq1})}.
\end{eqnarray*}

We have established that, apart from (\ref{knotleq}), the equations of
Table \ref{algtable} represent a system of equations equivalent to 
(\ref{genTrepeqcom}) and (\ref{genTrepeqadj}).

\np
We now consider the coefficients of $t^2$ and of $t^3$
in the equations (\ref{Artinrels}) for $\s_i$.
We claim that, given (\ref{commlaweq})--(\ref{1-1eq3}),
a necessary condition for the corresponding
equations to hold is
\begin{eqnarray}\label{knotleq}
T_{k,\a_l} &=& 0 \quad\ \mbox{if }\ k\ne l.
\end{eqnarray}
To see this, note that,
if $k\not\sim l$,
the coefficient of $t^2$ gives
$T_kT_l=T_lT_k$
which, applied to $x_{\a_l}$,
yields (\ref{knotleq}).  If $k\sim l$, note 
$$T_k\tau_lx_{\a_k}=
T_k((1-r^2)x_{\a_k}+rx_{\a_k+\a_l})=(r^4(1-r^2)+rT_{k,\a_k+\a_l})x_{\a_k}
=0$$ as $T_{k,\a_k+\a_l}=r^5-r^3$.  Now use the action of
$$T_l\tau_k T_l + \tau_lT_kT_l+T_lT_k\tau_l= T_k\tau_l T_k + \tau_k
T_lT_k+T_kT_l\tau_k $$
on $x_{\a_l}$.  
We see only the middle terms do not vanish because of the relation above
and so 
$$r^4T_{k,\a_l}\tau_lx_{\a_k}=T_{k,\a_l}T_{l,\a_k}\tau_kx_{\a_l}.$$ By
considering the coefficient of $x_{\a_k}$, which occurs only on the
left hand side, we see that (\ref{knotleq}) holds.

A consequence of this is that $T_iT_j=0$ if $i\neq j$.  Now all the
equations for the $t^2$ and $t^3$ coefficients are easily satisfied.
In the noncommuting case of $t^2$, the first
terms on either side are $0$ by the relation
above and the other terms are $0$ as $T_lT_k=0$.

It remains to establish that the matrices $\s_k$ are invertible.  To
prove this, 
we observe that the linear transformation $\s_k^2+(r^2-1)\s_k-r^2$ maps $V$
onto the
submodule spanned by $x_{\a_k}$ and that the image
of  $x_{\a_k}$ under $\s_k$ is $tr^4x_{\a_k}$.
In fact, the determinant of $\s_k$
equals $(-1)^c tr^{4+2c}$, where $c$ is the number of positive roots
$\b$ such that $(\a_k,\b) = -1$.
\end{proof}

\np
For a positive root $\b$, we write $\het(\b)$ to denote its height,
that is, the sum of its coefficients with respect to the $\a_i$.
Also, $\Supp(\b)$ is the set of $k\in\{1,\ldots,n\}$ such that
the coefficient of $\a_k$ in $\b$ is nonzero.

\begin{Cor}\label{lastline}
If the $T_{k,\b}$ satisfy the equations in Table \ref{algtable}, then
they obey the following rules, where $\het(\b)$ stands for the height
of $\b$.
\begin{enumerate}[(i)]
\item 
If $(\a_k,\b)=(\a_l,\b)=0$
and $(\a_k,\a_l)=-1$, then $T_{k,\b} = T_{l,\b}$.
\item If $(\a_k,\b) = 1$, then
$T_{k,\b} = r^{\het(\b)+1}(r^2-1)$.
\item The degree of
$T_{k,\beta}$ equals $3+{\rm ht}(\beta)$ whenever $k\in \Supp(\b)$.
\item $T_{k,\beta}$ is a
multiple of $r^2-1$ whenever $\beta\ne\a_k$.
\item $T_{k,\b}=0$  whenever $k\not\in \Supp(\b)$.
\end{enumerate}
\end{Cor}

\np
\begin{proof}
(i). The equations are necessary as they appeared under (\ref{00eq}).

\nl(ii).
Use induction on $\het(\b)$. If $\het(\b) = 2$, the equation coincides with
(\ref{11eq1}). If $\het(\b) >2$, then either (\ref{commlaweq})
or (\ref{01eq2}) applies.

\nl(iii) and (iv) are obvious.

\nl(v) follows from (\ref{knotleq}) by use of (\ref{commlaweq}) and
(\ref{1-1eq1}). Observe that, if $k\not\in\Supp(\b)$ and $(\a_l,\b) =1$
for some $l\sim k$, then $l\not\in\Supp(\b-\a_l)$.
\end{proof}

\np
The proposition enables us to describe
an algorithm computing the $T_{k,\beta}$, and which shows that there
is at most one solution.

\begin{Alg}\label{alg}
The polynomials $T_{k,\b}$ of Theorem \ref{adethm} can be computed 
as follows by using Table \ref{algtable}.
\begin{enumerate}[(i)]
\item
If $k\not\in\Supp(\b)$, then $T_{k,\b} = 0$.
Otherwise, proceed with the next steps.
\item
If ${\rm ht}(\beta)\le2$, equations
(\ref{00eq}) and (\ref{11eq1}),
that is, the first two lines of Table \ref{algtable},
determine
$T_{k,\beta}$.
\\
{From} now on, assume 
${\rm ht}(\beta)>2$.
We proceed by recursion, expressing $T_{k,\b}$
in $\Z[r]$-linear combinations of 
$T_{m,\gamma}$'s with ${\rm ht}(\gamma)<{\rm ht}(\b)$.
\item
Search for an $l\in\{1,\ldots,n\}$ such that $(\a_k,\a_l)=0$ (so $k$
and $l$ are nonadjacent in $M$) and $(\a_l,\b)=1$ (so
$\beta-\a_l\in\Phi$).
If such an $l$ exists, then (\ref{commlaweq})
expresses $T_{k,\beta} $ as a multiple of $T_{k,\b-\a_l}$.
\item
So, suppose there is no such $l$. Then search for an
$l\in\{1,\ldots,n\}$ such that $(\a_k,\a_l)=-1$ (so $k$
and $l$ are adjacent in $M$) and $(\a_l,\b)=1$ (so
$\b-\a_l\in\Phi$).
If such an $l$ exists, we cannot have $(\a_k,\b)=1$, for otherwise
$(\b,\a_k+\a_l)=2$ so $\b=\a_k+\a_l$, contradicting ${\rm ht}(\b)>3$.
Hence $(\a_k,\b)=0$ or $-1$; the identities
(\ref{01eq1}) and  (\ref{1-1eq1})
then express $T_{k,\beta}$ 
as a linear combination of $T_{k,\b-\a_l}$ and some $T_{l,\gamma}$
with ${\rm ht}(\gamma)<{\rm ht}(\b)$.
\item
Since there always exists an $l\in\{1,\ldots,n\}$ with $(\a_l,\b)=1$,
we are left with the case where $(\a_k,\b)=1$ and $(\a_m,\b)\ne1$ for
all $m\ne k$.  Now take $l$ adjacent to $k$.
Then $(\a_l,\b)=0$ or $-1$.  We argue that there is one $l$ for which
$(\a_l,\b)=0$ in which case (\ref{01eq2}) applies.  Indeed there is
some $l$ for which $(\a_l,\b-\a_k)=1$ as $\het(\b)\ge2$.  It cannot be
a root for which $(\a_k,\a_l)=0$ for then $1 =
(\a_l,\b-\a_k)=(\a_l,\b)\in\{-1,0\}$.  This means $(\a_k,\a_l)=-1$ and
so $(\a_l,\b)=0$.
\end{enumerate}
This ends the algorithm.
Observe that all lines of Table \ref{algtable} have been used.
\end{Alg}

\np
The algorithm computes a polynomial for each $k,\b$
based on Table \ref{algtable}, showing that there is at most one solution
to the set of equations. The next result shows that the computed polynomials
are indeed a solution.

\np
\begin{Prop}\label{uniquesol}
The equations of Table \ref{algtable} have a unique solution.
\end{Prop}

\np
\begin{proof}
We need to show that the matrices $\s_i$ satisfy the relations in
Proposition \ref{finalTiRelations}.  It is clear from (\ref{11eq2})
and (\ref{11eq1}) that they do for $\b$ of height $1$ or $2$.  We use
induction on $\het(\b)$, the height of $\b$, and assume $\het(\b)\ge3$.

We first check (\ref{commlaweq}).  If this applies, the value
$T_{k,\b}$ was determined in step (iii) of the algorithm, and we are
really checking the value did not depend on the choice of $l$.
Suppose that there are two choices $l$ and $l'$ for which
$(\a_l,\b)=(\a_{l'},\b)=1$ and $(\a_l,\a_k)=(\a_{l'},\b)=0$.  Then by
our definition $T_{k,\b}=rT_{k,\b-\a_l}$ and we must show that
$T_{k,\b}=rT_{k,\b-\a_{l'}}$.  Suppose first $(\a_l,\a_{l'})=0$.  Then
$(\b-\a_l,\a_{l'})=1-0=1$ and $(\b-\a_{l'},\a_l)=1$.  This means
$\b-\a_l-\a_{l'}$ is also a root.  Now apply (\ref{commlaweq}) and the
induction to see $T_{k,\b-\a_l}=rT_{k,\b-\a_l-\a_{l'}}$ and
$T_{k,\b-\a_{l'}}=rT_{\b-\a_l-\a_{l'}}$, and so they are the same.

We have dealt with the other case before in which $(\a_l,\a_{l'})=-1$
and seen this only occurs when $\b$ has height $2$.

We now consider the case in which there are $l$, $l'$ with
$(\b,\a_l)=(\b,\a_{l'})=1$, $(\a_k,\a_l)=-1$ and $(\a_k,\a_{l'})=0$.
Here $T_{k,\b}$ is determined by step (iii) of the algorithm,
$T_{k,\b}=rT_{k,\b-\a_{l'}}$.  We
must have $(\a_l,\a_{l'})=0$, for if $(\a_l,\a_{l'})=-1$, then we would again
be in the height $2$ case as $(\b-\a_l,\a_{l'})=2$.  We must show for
(\ref{01eq1}) that if $ (\b,\a_l)=0$ then
$$rT_{k,\b-\a_{l'}}=T_{l,\b-\a_k-\a_l}+(r-r^{-1})T_{k,\b-\a_l}.$$
Observe that $(\b-\a_{l'},\a_l)=1$ and $(\a_k,\a_l)=-1$ and so, by
(\ref{01eq1}),
$$rT_{k,\b-\a_{l'}}=rT_{l,\b-\a_{l'}-\a_l-a_k}
+r(r-r^{-1})T_{k,\b-\a_{l'}-\a_l}.$$ Now, as $(\a_l,\a_{l'})=0$, we can
use (\ref{commlaweq}) to obtain
$T_{l,\b-\a_k-\a_l}=rT_{l,\b-\a_l-\a_k-\a_{l'}}$ and
$T_{k,\b-\a_l}=rT_{k,\b-\a_l-\a_{l'}}$, and so the equations are satisfied.
In order to satisfy (\ref{1-1eq1}) when $(\b,\a_l)=-1$, we need to show
$$rT_{k,\b-\a_{l'}}=rr^{-1}T_{l,\b-\a_l}+r(r-r^{-1})T_{k,\b-\a_l}.$$
Again express these terms using (\ref{commlaweq}) subtracting $\a_{l'}$ in 
each of the expressions to get 
equality.    

We have now checked all equations except (\ref{01eq2}) keeping in mind
that (\ref{00eq}) will be satisfied if the others are.  Now
(\ref{01eq2}) follows from
Corollary \ref{lastline}(ii) as $(\b-\a_k,\a_l)=1$. 
\end{proof}

\np 
We are now ready to prove the first part of Theorem \ref{adethm}.

\begin{Cor}\label{proofthmpart1}
The solution $T_{k,\b}$ described by Proposition \ref{uniquesol}
determines a representation of $B$ on $V$ as claimed in Theorem
\ref{adethm}.
\end{Cor}

\np
\begin{proof}
Combine Propositions \ref{finalTiRelations} and \ref{uniquesol}.
\end{proof}

\begin{Cor}\label{algcomb}
The solution $T_{k,\b}$ of Proposition \ref{uniquesol} is computable via
exponents $a_{k,\b}$, $c_{k,\b}$, $d_{k,\b}$ as follows.
$T_{k,\b} = 0$ if $\a_k\not\in{\rm Supp}(\b)$. Otherwise,
\begin{eqnarray}\label{Tcomdf}
\frac{T_{k,\b}}{r^{\het(\b)+1}(r^2-1)}
&=& 
\left\{
\begin{tabular}{lr}
$1$ &\mbox{ if }\ $\phantom{-}(\a_k,\b)=1$\\
$(1-r^{-a_{k,\b}})$ & \mbox{ if }\ $\phantom{-}(\a_k,\b)=0$\\
$(1-r^{-c_{k,\b}})(1-r^{-d_{k,\b}})$ & \mbox{ if }\ $(\a_k,\b)=-1$
\end{tabular}
\right.
\end{eqnarray}
where $a_{k,\b}$, $c_{k,\b}$, $d_{k,\b}$
are $0$ if $k\not\in\Supp(\b)$; otherwise,
\begin{eqnarray*}
a_{k,\b} &=& a_{k,\b-\a_l} \phantom{24567890moremm}\,\, \mbox{ if } (\a_l,\b)=1
\mbox{ and } k\not\sim l\\
a_{k,\b} &=& a_{l,\b-\a_l-\a_k}+2 \phantom{2344567890}\,\, \mbox{ if } 
(\a_l,\b)=1
\mbox{ and } k\sim l\\
\{c_{k,\b},d_{k,\b}\} &=& \{c_{k,\b-\a_l},
d_{k,\b-\a_l}\} \phantom{mremm}\,\, \mbox{ if } (\a_l,\b)=1
\mbox{ and } k\not\sim l\\
\{c_{k,\b},d_{k,\b}\} &=& \{a_{l,\b-\a_l},c_{l,\b-\a_l}+2\} 
 \phantom{12345} \mbox{ if } 
(\a_l,\b)=1,\  d_{l,\b-\a_l} =
a_{k,\b-\a_l},\\
\nonumber&& \phantom{alotttttt23456789moremm}  \mbox{ and } \ k\sim l\\  
\{c_{k,\b},d_{k,\b}\} &=& \{a_{l,\b},a_{m,\b}\} \phantom{13456789more}\,\,
 \mbox{ if } 
(\a_l,\b)=0,\ (\a_m,\b)=0, \\
\nonumber&& \phantom{alotttttt23456788moremm}  \mbox{ and } \ l\sim k\sim m\ne l.
\end{eqnarray*}
\end{Cor}

\np
\begin{proof}
The proof is similar to that of Proposition
\ref{finalTiRelations}.
\end{proof}

\np
\begin{Example}
The $A_n$ case.
Then $c_{k,\b}=d_{k,\b} = 0$ and $a_{k,\b} = 2$ if $k\in\Supp(\b)$.
Note that the last three lines of the corollary do not occur.
Our representation can be obtained from
the Lawrence Krammer representation
as described in \cite{Krammer}
by a diagonal transformation with respect to the basis
$x_\b$ $(\b\in \Phi^+)$,
and by replacing $q$ by $r^2$. As a result, the involutory
automorphism of the diagram $A_n$ can be realized as a linear
transformation leaving invariant the basis (compare with Remark 5.1 of
\cite{Krammer}).
\end{Example}

\begin{Example} The $D_n$ case.
For a given root in $D_n$, let $l_1$ be the number of coefficients $1$
in the expression of $\b$ as a linear combination of the $\a_i$, and
let $l_2$ be the number of coefficients $2$.  These are the only
nonzero coefficients which can occur for $D_n$.  In the case
$(\a_k,\b)=-1$, the exponent $c_{k,\b}$ is $2$ and the exponent
$d_{k,\b}$ is $2l_2+2$.  Assume now that $(\a_k,\b)=0$ and
$k\in\Supp(\b)$. If $\a_k$ has coefficient 2 in $\b$ or $k$ is the end
node of a short branch of the Coxeter diagram, then $a_{k,\b} = 4$; if
$k$ is the end node of the long branch (possibly after removing nodes
with zero coefficients), then $a_{k,\b} = 2l_2+2$; otherwise $a_{k,\b}
= 2$.  It is straightforward to check that the relations of Table
\ref{algtable} all hold.
\end{Example}

\section{Faithfulness of the representation}
\label{krammer}

We now combine the representation of Section \ref{rep} with the
root system knowledge of Section \ref{basics}.  Our arguments are
straightforward generalizations of Krammer's method, but we give
details anyway for the reader's convenience.

Recall that $V$ is the free module over $\Z[t^{\pm1},r^{\pm1}]$ generated
by $x_\b$ for $\b$ ranging over the positive roots.  In Corollary
\ref{proofthmpart1} we established the first part of Theorem \ref{adethm}. In
this section we prove the second part.  To this end, we specialize
$r$ to a real number $r_0$ with $0<r_0<1$ in
$V\otimes\R$ to obtain $V_1$,
the free module over $\R[t,t^{-1}]$ generated by the $x_\b$.
We also keep the Coxeter matrix $M$ to be one of
$A_n$
$(n\ge1)$, $D_n$ $(n\ge4)$, $E_6$, $E_7$, or $E_8$.

Note that $0<r_0<1$ implies that the constant term
of each of the entries of the matrices
$\s_i$ is a nonnegative real number.  This will be the same for any
product of $\s_i$, and so for any element of the monoid $B^+$ they
generate.  Therefore, in its linear action on $V_1$, the monoid $B^+$
preserves
\begin{eqnarray}
\label{Udef}
U &=& \oplus_{\b\in \Phi^+}  (\R_{\ge0}\oplus t\R[t])x_\b.
\end{eqnarray}
For $A\subseteq \Phi^+$ set
$$U_A =\left\{\left.\sum_{\b\in\Phi^+} k_\b 
x_\b\in U \,\right|\, k_\b\in t\R[t]\Leftrightarrow\b\in A\right\}.$$
Then, obviously, $U$ is the disjoint union of the $U_A$.

\begin{lemma}
For $x\in B^+$ and $A\subseteq\Phi^+$, there is a unique $A'\subseteq\Phi^+$
such that $xU_A\subseteq U_{A'}$.
\end{lemma}

\begin{proof}
For a given subset $A$ of $\Phi^+$, the elements of $U_A$ are the
vectors in $U$ for which the support mod $t$ is exactly
$\Phi^+\setminus A$.  In
particular, an element $u=\sum_{\b\in\Phi^+} (k_\b+tp_\b)x_\b$ of $U$,
with $k_\b\in\R_{\ge0}$ and $p_\b\in\R[t]$, is in $U_A$ if and only if
$k_\b=0$ for $\b\in A$ and $k_\b\neq 0$ for $\b\in\Phi^+\setminus A$.
As all matrix entries of an element $x$ of $B^+$ are nonnegative mod
$t$, the image by $x$ acting on two nonzero elements of $U_A$ will
have exactly the same support mod $t$.  If this is $\Phi^+\setminus
A'$, the images of nonzero vectors of $U_A$ are all in $U_{A'}$.
\end{proof}

\np The assignment $(x,A)\mapsto A'$ where $A'$ is the unique subset
of $\Phi^+$ such that $U_{A'}$ contains $xU_A$, defines an action of
$B^+$ on ${\mathcal P}(\Phi^+)$; we write $x\sact A$ for $A'$.
Observe that $A\subseteq D $ implies that $ x\sact A\subseteq x\sact
D$.

\begin{lemma}\label{sactclosed}
\label{s_istacAlm}
The action $\sact$ preserves $\mathcal C$.  
It can be explicitly described for $s_k$ as follows,
where $k\in\{1,\ldots,n\}$ and $A\in{\mathcal C}$.
\begin{eqnarray*}
s_k\sact A= \{\a_k\}&\cup&
\left\{\b\in\Phi^+\left|
\begin{tabular}{lr}
${\b-\a_k}\in A$ & \mbox{ if }\quad $\phantom{-}(\a_k,\b)=1$,\\
${\b}\in A$ & \mbox{ if }\quad $\phantom{-}(\a_k,\b)=0$,\\
$ {\b} ,{\b+\a_k}\in A$ &\phantom{-}
 \mbox{ if }\quad $(\a_k,\b)=-1$\phantom{,}\\
\end{tabular}
\right.\right\}
\end{eqnarray*}
In particular, $\a_k\in s_k\sact A \subseteq \{\a_k\}\cup r_k(A)$.
\end{lemma}

\np
\begin{proof}
For the proof of the first statement, it suffices to consider $x=s_k$
as $B^+$ is generated by these elements. 

As for the description of $s_k\sact A$, only the action of 
$\tau_k$ on $u=\sum_{\b\in\Phi^+}k_\b x_\b\in U$ 
with $k_\b\in \R_{\ge0}$ is relevant.
A computation shows
\begin{eqnarray*}
\tau_ku& = &
\sum_{(\c,\a_k) = -1}k_\c((1-r_0^2)x_\c + r_0x_{\c+a_k})
+
\sum_{(\c,\a_k) = 0}k_\c x_{\c}
+
\sum_{(\c,\a_k) = 1}k_\c r_0 x_{\c-a_k}
\\
&=&
\sum_{(\b,\a_k) = 1}k_{\b-\a_k} r_0x_{\b}
+
\sum_{(\b,\a_k) = 0}k_\b x_{\b}
+ \sum_{(\b,\a_k) =- 1} (k_{\b+\a_k}r_0+k_{\b}(1-r_0^2)) x_{\b}.
\end{eqnarray*}
The set $s_k\sact A$ is the set of positive
roots for which $x_\b$ has coefficient 
$0$ in $\tau_ku$ for any element $u$ in $U_A$.
The description of $s_k\sact A$ follows 
directly from this formula.
For instance, for $\beta\in\Phi^+$ with $(\b,\a_k) =- 1$
to belong to $s_k\sact A$, we need to have
$k_{\b+\a_k} r_0 + k_\b(1-r_0^2) = 0$, which is equivalent to
$k_{\b+\a_k} =k_\b = 0$, whence $\b+\a_k,\b\in A$. 

It remains to show that $s_k\sact A$ is closed.  So suppose that $\b$
and $\c$ are in $s_k\sact A$ and that $\b+\c$ is in $\Phi^+$. 
Assume $\c=\a_k$.  We always have $\a_k $ in $s_k\sact
A$.  As $\b+\a_k\in\Phi^+$, the inner product 
$(\a_k,\b)$ equals $-1$. By the above, this implies
that both $\b$ and $\b+\a_k$ are in $A$.
But then $\b+\a_k\in\Phi^+$ satisfies $(\a_k,\b+\a_k) = 1$ and
$(\b+\a_k)-\a_k\in A$, so $\b+\a_k\in s_k\sact A$.

{From} now on, we assume that neither $\b$ nor $\c$ are equal to $\a_k$.
Suppose that both $\b$ and $\c$ are orthogonal to
$\a_k$.  We saw above that being in
$s_k\sact A$ means that both $\b$ and $\c $ are in $A$ and because $A$
is closed, $\b+\c$ is also in $A$.  But then $\b+\c $, being
orthogonal to $\a_k$, also belongs to $s_k\sact A$.

The case remains where at least one of $\b$ and $\c$ is not orthogonal
to $\a_k$.  Suppose first that $(\a_k,\b)=-1$.  As $\b\in s_k\sact A$,
by the above, both $\b$ and $\b+\a_k$ are in $A$.  If $\c$ is
orthogonal to $\a_k$ we know from above and from $\c \in s_k\sact A$
that $\c $ is in $A$.  Now, as $\b+\a_k,\c\in A$ and $A$ is closed,
also $\c+\b+\a_k\in A$.  As $\b,\c\in A$ and $A$ is closed, also
$\b+\c$ is in $A$.  Now $(\b+\c,\a_k)=-1$ and so by the above $\b+\c$
is in $A$.  We still need to consider the other possibilities for
$(\c,\a_k)$.  As $\b+\c$ is a root, $(\c,\a_k) \neq -1$.  Now
$(\b+\c,\a_k)=0$ and we need only show that $\b+\c\in A$.  But this
follows as $A$ is closed and $\b+\a_k,\c-\a_k\in A$.

The only case remaining is $(\b,\a_k)=1$ and $(\c,\a_k)\in\{0, 1\}$.
However, the latter inner product cannot be 1, for otherwise 
$(\b+\c,\a_k)=2$, contradicting the fact that $\b+\c$ is a
positive root.  This means $(\c,\a_k)=0$ and as $\c\in s_k\sact A$, we find
$\c\in A$.  As $(\b,\a_k)=1$ and $\b\in s_k\sact A$, we have
$\b-\a_k\in A$.  Since $(\b+\c,\a_k)=1$, the vector $\b+\c-\a_k$ is a
positive root.  As both $\b-\a_k$ and $\c$ are in $A$ and $A$ is closed, the
root $\b-\a_k+\c$ belongs to $A$.  Now as $(\b+\c,\a_k)=1$ and
$\b+\c-\a_k\in A$, we conclude $\b+\c\in s_k \sact A$.
\end{proof}

\np
\begin{Lm}\label{4.4lm}
For $w\in W$ and $i\in\{1,\ldots,n\}$ satisfying
$l(r_iw)<l(w)$, and for each closed subset $A$ of $\Phi^+$,
we have
$\Phi_w\subseteq \{\a_i\}\cup r_i(A)$ if and only if
$w\le b^{-1}(L(s_ig(A)))$.
\end{Lm}

\np
\begin{proof}
Since $l(r_iw)<l(w)$, the subset $\Phi_w$ of $\Phi^+$ 
coincides with $\{\a_i\}\cup r_i(\Phi_{r_iw})$.
Hence $\Phi_w\subseteq \{\a_i\}\cup r_i(A)$ if and only if
$\Phi_{r_iw}\subseteq A$, which, by definition of $g$,
is equivalent to $b(r_iw)\le g(A)$. By Proposition
\ref{Artinmonprops}(iii), this is the same as 
$s_ib(r_iw)\le s_ig(A)$, while, since the left hand side equals
$b(w)$, this in turn amounts to $b(w)\le L(s_ig(A))$.
Hence the lemma.
\end{proof}

\np
\begin{Lm}\label{4.5lm}
Suppose that the subsets  $A$ and $E$ of $\Phi^+$ are
closed and, for some $i\in\{1,\ldots,n\}$,
satisfy 
$\{\a_i\}\subseteq E \subseteq  \{\a_i\} \cup r_i(A)$.
Then $E  \subseteq  s_i\sact A$.
\end{Lm}

\np
\begin{proof}
Let $\b\in E$.  We show that $\b\in s_i\sact A$.
We distinguish cases according to $(\a_i,\b)$.
If $(\a_i,\b)=2$, then $\b=\a_i\in s_i\sact A$ by Lemma \ref{sactclosed}.

\nl
If $(\a_i,\b)=1$, then $\b=r_i(\b-\a_i)$ 
with $\b-\a_i\in A$. By Lemma \ref{s_istacAlm}, this implies
$\b\in s_i\sact A$.

\nl
If $(\a_i,\b)=0$, then $\b=r_i(\b)$ 
with $\b \in A$. By Lemma \ref{s_istacAlm}, this implies
$\b\in s_i\sact A$.

\nl
Finally, suppose $(\a_i,\b) = -1$. Then $\b = r_i(\b+\a_i)$ with $\b+\a_i\in A$.
Moreover, since $\b,\a_i\in E$ and $E$ is closed, we have $\b+\a_i\in E$,
which, by analysis of the case of inner product $1$ above,
gives $\b+\a_i\in s_i\sact A$. In particular, 
$\b+\a_i\in r_i(A)$, so $\b+\a_i=r_i(\b)$ with $\b\in A$.
Since $\b,\b+\a_i\in A$, Lemma \ref{s_istacAlm} implies
$\b\in s_i\sact A$.
Hence the lemma.
\end{proof}

\np
\begin{Prop}\label{equivariantprop}
The map $g:{\mathcal C} \to\Omega$ is $B^+$ equivariant.
That is, for all $x\in B^+$ and $A\in {\mathcal C}$,
we have
$$g(x\sact A) = L(x g(A)).$$
\end{Prop}

\np
\begin{proof}  It suffices to prove the assertion
for $x=s_i$ with $1\le i \le n$.  Write
$w= b^{-1}g(s_i\sact A)$,
so $b(w) = g(s_i\sact A)$, and $\Phi_{w}$ is the maximal subset of
$s_i\sact A$ belonging to $\mathcal D$.
Recall from Lemma \ref{sactclosed} that $\a_i\in s_i \sact A$.
It implies $\Phi_{r_i}\subseteq s_i \sact A$,
whence $\Phi_{r_i}\subseteq \Phi_w$,
so $r_i\le w$. In other words, $l(r_iw)<l(w)$.
Since $s_i\sact A\subseteq \{\a_i\}\cup r_i(A)$,
we obtain 
$\Phi_{w}\subseteq \{\a_i\}\cup r_i(A)$.
Put $w' = b^{-1}(L(s_ig(A)))$.
By Lemma \ref{4.4lm}, $\Phi_{w'}$ is the maximal
element of $\mathcal D$ contained in
$\{\a_i\}\cup r_i(A)$ and, by Lemma \ref{4.5lm}, so is $\Phi_w$.
Therefore, by Lemma \ref{Phiprops}(v), $w=w'$, proving
$g(s_i\sact A) = b(w) = L(s_ig(A))$.
\end{proof}

\np
For $x\in\Omega$,
write 
\begin{eqnarray}\label{Cdef}
C_x &=& \bigcup_{A\in {\mathcal C}, g(A)=x} U_A.
\end{eqnarray}

\np
\begin{Prop}\label{CxProp}
The subsets $C_x$ $(x\in\Omega)$ satisfy the following three properties
for each $x,y\in\Omega$.
\begin{enumerate}[(i)]
\item $C_x\ne\emptyset$.
\item $C_x\cap C_y=\emptyset$  if $x\ne y$.
\item $xC_y\subset C_{L(xy)}$.
\end{enumerate}
\end{Prop}

\np
\begin{proof}
(i).
Clearly, $\emptyset\ne U_{\Phi_{b^{-1}(x)}}\subseteq C_x$, so 
$C_x$ is nonempty.

\nl (ii).
This follows immediately from the definition of $C_x$.

\nl(iii).
Given $x,y\in\Omega$,
let $A\in{\mathcal C}$ be such that $y=g(A)$.
Then, by respectively, the definition of $\sact$,
the definition of $C_x$, and Proposition \ref{equivariantprop}(vi),
$$xU_A\subseteq U_{x\sact A}\subseteq C_{g(x\sact A)} =
C_{L(xy)},$$
whence
$xC_y\subseteq C_{L(xy)}$.
\end{proof}

\np
In fact, (iii) also holds for each $x\in B^+$, as follows
from the following argument based on induction with respect to $l(x)$.
If $l(x)>1$ then
there exist $i\in\{1,\ldots,n\}$ 
and $u\in \Omega$ such that $x=s_iu$ and $l(x) = 1+l(u)$.
Then, by  the induction hypothesis, (iii) of the proposition,
and Proposition \ref{Artinmonprops},
$$
xC_y =
s_iuC_y \subseteq s_iC_{L(uy)}\subseteq C_{L(s_iL(uy))}
= C_{L(s_iuy)}= C_{L(xy)}.
$$

\np
\begin{Prop}\label{Cxfaithful}
Let $B^+$ act on a set $U$ in such a way that each element acts injectively.
Suppose we are given subsets $C_x$ of $U$ for $x\in \Omega$
satisfying properties (i), (ii), and (iii) of Proposition \ref{CxProp}.
Then the action of $B^+$ on $U$ is faithful.
\end{Prop}

\np
\begin{proof}
(This is the proof appearing in \cite{Krammer}.)
Suppose that the elements $x$
and $y$ of $B^+$ act identically on $U$.  If $l(x)+l(y)=0$, then $x$
and $y$ are both the identity and there is nothing to prove.
Suppose therefore, that $l(x)+l(y)>0$.
Pick $u\in C_1$. Then $xu\in xC_1\cap yC_1\subseteq C_{L(x)}\cap C_{L(y)}$,
which implies by Proposition \ref{CxProp} that 
$z = {L(x)}= {L(y)}$ for some nontrivial $z\in\Omega$. This means that
there are $x'$, $y'$ in $B^+$ such that
$x=zx'$ and $y=zy'$. But then,
as $z$ acts injectively, $x'$ and $y'$ act identically on $U$, whereas
$l(x')+l(y')=l(x)+l(y)-2l(z)$, so we can finish
by induction on $l(x)+l(y)$.
\end{proof}

\np
{\bf Proofs of Theorems \ref{mainthm} and \ref{adethm}.}
Propositions \ref{CxProp} and \ref{Cxfaithful}
with $U$ as in (\ref{Udef}) and
$C_x$ as in (\ref{Cdef}),
together with Corollary \ref{proofthmpart1}, give a proof of
Theorem \ref{adethm}.  As for Theorem \ref{mainthm}, suppose that $M$
is of finite type.  If $M$ is the disjoint union of diagrams $M'$ and
$M''$, then $B$ is the direct product of the Artin groups $B'$, $B''$
corresponding to $M'$, $M''$, respectively, and so the direct sum of
faithful linear representations of $B'$ and $B''$ would be a faithful
linear representation for $B$.  Hence, a proof of Theorem
\ref{mainthm} in the case where $M$ is finite and irreducible
suffices.

By \cite{Crisp}, every Artin group $B$ of finite type $M$ 
such that $M$ has a multiple bond occurs as a subgroup of an Artin group
of finite type without multiple bonds.
Therefore, a proof of Theorem \ref{mainthm} for finite irreducible types
without multiple bonds, that is, for types $A$, $D$, $E$, suffices,
and this is dealt with by Theorem \ref{adethm}.
This ends the proof of the theorems in Section \ref{intro}.

\section{Epilog}
\label{epilog}
As stated before, the Artin groups whose types are spherical
irreducible Coxeter matrices with multiple bonds occur as subgroups of
Artin groups of finite types without multiple bonds.  They occur as
fixed subgroups of an automorphism group $H$ of $B$ permuting the
vertices of $M$.  The natural generators of this subgroup are the
elements $\prod_{k\in E}s_k$ of $B$ for $E$ running over the $H$
orbits on the vertex set of $M$.  It is obvious that these subgroups
satisfy the Artin group relations, but it is harder to establish that
every relation they satisfy is a consequence of these. It may be of
interest to know whether the latter can also be proved by applying
Krammer's methods to the representation of the $H$ fixed
subgroup of $B$ on the centralizer in $V$ of $H$ with respect to a
suitable action of $H$ on $V$.

We are also able to recover Theorem~$6.1$ from \cite{Krammer}.  For
this we need the {\it Charney} length function, $l_\Omega$ on $B$.  It
assigns to $x\in B$ the smallest natural number $k$ such that there
are elements $x_1,x_2,\dots,x_k$ in $\Omega \cup \Omega^{-1}$ for
which $x=x_1x_2\cdots x_k$.

\begin{Thm}\label{Urmk}
Let $B$, $R$, $V$ be as in Theorem
\ref{adethm} and write $\rho$ for the linear representation
$B\to {\rm GL}(V)$.
For
$x\in B$, consider the Laurent expansion of $\rho(x)$ with respect to $t$:
$$\rho(x)
=
\sum_{i=k}^h A_it^i, \quad \ A_k\ne0,\ \ A_h\ne 0,$$
where $A_i$ is a matrix whose entries are in $\Z[r^{\pm1}]$.
\begin{enumerate}[(i)]
\item Then $l_\Omega(x)=\max(h-k,h,-k)$.
\item 
If in addition $x\in B^+\setminus b(w_0)B^+$, then $k=0$ and $h=l_\Omega(x)$.
Here $w_0$ is the longest word in the Coxeter group $W$ corresponding to $B$.
\end{enumerate}
\end{Thm}

\np
\begin{proof}
The proof is as in \cite{Krammer} and so we do not include it.
The use of Lemmas~$3.1$ and $3.2$ in \cite{Krammer} is
replaced by the following corresponding results for $\rho$.

There is a linear transformation $U\in {\rm GL}(V)$ whose
matrix with respect to $\{x_\b\}_\b$ have entries in $\Z[r^{\pm1}]$
such that $\s_k U \widehat \s_k = U$ for each $k\in\{1,\ldots,n\}$,
where $\widehat \s_k$ is the matrix $\s_k$ with $t$ and $r$ 
replaced by $t^{-1}$ and $r^{-1}$, respectively.
The matrix $U$ is determined by
the following rules involving an index $k\in\{1,\ldots,n\}$ such that
$(\a_k,\b) = 1$.
\begin{eqnarray*}
U_{\c,\b} &=&
\left\{
\begin{tabular}{lcr}
$0$ &\mbox{ if }&\quad $\phantom{-}\c\not\le\b$,\\
$1$ &\mbox{ if }&\quad $\phantom{-}\c=\b$,\\
$\widehat T_{k,\b} r^4$ & \mbox{ if }&\quad $\c=\a_k\le\b$,\\
$ U_{\c-\a_k,\b-\a_k}$ &
 \mbox{ if }&\quad $\c\le\b$ and  $(\a_k,\c)=1$,\\
$r^{-1} U_{\c,\b-\a_k}$ &
 \mbox{ if }&\quad $\c\le\b$ and  $(\a_k,\c)=0$,\\
$ U_{\c+\a_k,\b-\a_k}+(r^{-1}-r)U_{\c,\b-\a_k}$ &
 \mbox{ if }&\quad $\c\le\b$ and $(\a_k,\c)=-1$.\\
\end{tabular}
\right.
\end{eqnarray*}
This matrix replaces the matrix $T(q)$ in Lemma~$3.1$ in 
\cite{Krammer}. 

In the representation of Theorem \ref{adethm}, $\rho(b(w_0))$ is the
multiple of the permutation matrix $\pi$ by the scalar $tr^{e+3}$.  Here 
$\pi$ permutes $\{x_\b\}_\b$ according to the action of $-w_0$
on $\Phi^+$ and $e$ is the number of positive roots that are not
orthogonal to a given root.  In particular $e+3 = 2(n+1)$ 
for $A_n$, $4(n-1)$ for
$D_n$, $24$ for $E_6$, $36$ for $E_7$, and $60$ for $E_8$.
Note that this is 
in accordance with the theorem we are proving as $k=h=1$.
The matrix $\rho(b(w_0))$
replaces the matrix 
of Lemma~$3.2$ in \cite{Krammer}.  
\end{proof}

\np
Just as in \cite{Krammer}, this leads to a different proof that $\rho$ is
faithful.  Indeed, if $x$ is in the kernel, we see $h=k=0$ and so
$l_\Omega(x)= 0$, establishing that $x$ is the unit element of $B$.

\enddocument